\documentclass[11pt]{amsart}
\usepackage{graphicx}
\usepackage{graphics}
\usepackage{amsmath}
\usepackage{amscd}
\usepackage{latexsym}
\usepackage{hyperref}

\begin{document}

\textwidth 6.2in
\textheight 7.6in
\evensidemargin .75in
\oddsidemargin.75in

\newtheorem{Thm}{Theorem}
\newtheorem{Lem}[Thm]{Lemma}
\newtheorem{Cor}[Thm]{Corollary}
\newtheorem{Prop}[Thm]{Proposition}
\newtheorem{Rm}{Remark}

\def\a{{\mathbb a}}
\def\C{{\mathbb C}}
\def\A{{\mathbb A}}
\def\B{{\mathbb B}}
\def\D{{\mathbb D}}
\def\E{{\mathbb E}}
\def\R{{\mathbb R}}
\def\P{{\mathbb P}}
\def\S{{\mathbb S}}
\def\Z{{\mathbb Z}}
\def\O{{\mathbb O}}
\def\H{{\mathbb H}}
\def\V{{\mathbb V}}
\def\Q{{\mathbb Q}}
\def\Cn{${\mathcal C}_n$}
\def\CM{\mathcal M}
\def\CG{\mathcal G}
\def\CH{\mathcal H}
\def\CT{\mathcal T}
\def\CF{\mathcal F}
\def\CA{\mathcal A}
\def\CB{\mathcal B}
\def\CD{\mathcal D}
\def\CP{\mathcal P}
\def\CS{\mathcal S}
\def\CZ{\mathcal Z}
\def\CE{\mathcal E}
\def\CL{\mathcal L}
\def\CV{\mathcal V}
\def\CW{\mathcal W}
\def\IC{\mathbb C}
\def\IF{\mathbb F}
\def\IK{\mathcal K}
\def\IL{\mathcal L}
\def\IP{\bf P}
\def\IR{\mathbb R}
\def\IZ{\mathbb Z}

\title{The Akhmedov-Park exotic  ${\C\P}^{2}\#3\bar{\C\P}^2$}
\author{Selman Akbulut}
\thanks{The author is partially supported by NSF grant DMS 0905917}
\address{Department  of Mathematics, Michigan State University,  MI, 48824}
\email{akbulut@math.msu.edu }
\date{\today}
\begin{abstract} 
Here we draw a handlebody picture for the exotic ${\C\P}^{2}\#3\bar{\C\P}^2$ constructed by Akhmedov and Park.

\end{abstract}

\date{}
\maketitle

\setcounter{section}{-1}

\vspace{-.35in}

\section{Introduction}

Akhmedov-Park manifold $M$ is a symplectic $4$-manifold which is an exotic copy of ${\C\P}^{2}\#3\bar{\C\P}^2$ (\cite{ap}, also see related \cite{ak}). It is  obtained from two codimension zero pieces  which are glued along their common bondaries: $$M=\tilde{E_{0}}\smile_{\partial} \tilde{E_{2}}$$ The two pieces are constructed as follows:  We start with the product of a genus $2$ surface and the torus $E = \Sigma_{2}\times T^{2}$. Let $<a_1,b_1,a_2,b_2>$ and $<C,D >$ be the standard circles generating the first homology of $\Sigma_{2}$ and $T^{2}$ respectively (the cores of the $1$-handles). Then $\tilde{E}$ is obtained from $E$ by doing ``Luttinger surgeries'' to the four subtori $(a_1\times C, a_1)$, $(b_1\times C, b_1)$,   $(a_2 \times C, C)$  $(a_2 \times D, D)$ (see Section 1 for Luttinger surgery).  Then
$$\tilde{E_{0}}=\tilde{E}- \Sigma_{2}\times D^2$$

To built the other piece, let $K\subset S^3$ be the trefoil knot, and $S^{3}_{0}(K)$ be the $3$-manifold obtained doing $0$-surgery to $S^3$ along $K$. Being a fibered knot, $K$ induces a fibration $T^2\to S^{3}_{0}(K)\to S^1$ and the fibration
$$T^2\to S^{3}_{0}(K)\times S^1\to T^2$$
Let $T^{2}_1$ and $T^{2}_2$ be the vertical (fiber) and the horizontal (section) tori of this fibration, intersecting at one point $p$ . Smoothing $T_{1}^{2}\cup T^{2}_{2}$ at $p$ gives an imbedded genus $2$ surface with self intersection $2$, hence by blowing up the total space twice (at points on this surface) we get a genus $2$ surface   $\Sigma_{2} \subset (S^{3}_{0}(K)\times S^1) \# 2\bar{\C\P}^2$ with trivial normal bundle. Then we define: $$\tilde{E_{2}}=(S^{3}_{0}(K)\times S^1) \# 2 \bar{\C\P}^2 - \Sigma_{2} \times D^2$$


\noindent In short  $M= \tilde{E}\;\sharp_{\;\Sigma_{2}} (S^{3}_{0}(K)\times S^1) \# 2 \bar{\C\P}^2 $ ($\sharp_{\;\Sigma_{2}}$ denotes fiber sum along $\Sigma_{2}$).  Alternatively, $M$ can be built by using   $Sym^{2} (\Sigma_{3})$ 
\cite{fps}, which is equivalent to this construction, since 
$Sym^{2} (\Sigma_{3})=
  \;(\Sigma_{2} \times T^{2}) \# \bar{\C\P}^2 \sharp_{\;\Sigma_{2}}  (T^4 \# \bar{\C\P}^2 )$  and by Luttinger surgeries this can be transform to $M$ defined above. Also \cite{bk} gives another  exotic ${\C\P}^{2}\#3\bar{\C\P}^2$ which turns out to be a version of this $M$ (we thank Anar Akhmedov for explaining these equivalences).  

{\Rm Chronologically, first  Fintushel-Stern had the idea of building exotic  ${\C\P}^{2}\#3\bar{\C\P}^2$ from $Sym^{2} (\Sigma_{3})$ by Luttinger surgeries, but they couldn't get their manifold to be simply connected. Then Akhmedov-Park came out with this exotic ${\C\P}^{2}\#3\bar{\C\P}^2$ \cite{ap} (subject of this paper). Later in \cite{fps} Fintushel-Park-Stern fixed the fundamental group problem in their  $Sym^{2} (\Sigma_{3})$ approach, thereby getting another exotic ${\C\P}^{2}\#3\bar{\C\P}^2$ and in \cite{bk}  Baldridge-Kirk came out with their version. In retrospect, they all are related.}

\section{Luttinger surgery}

 Let $T^{2}\subset X$ be a smooth $4$-manifold with an imbedded subtorus which has trivial normal bundle $\nu (T^{2})\approx T^{2}\times B^2$.  Let $\varphi_p$ $(p\geq 0)$ be the self-diffeomorphism of $ T^{3}=T^{2}\times \partial B^{2}$ given by the matrix in terms of the automorphism of its standard homology generators $(a,b,c)$. 
\begin{equation*}
\left(
\begin{array}{ccc}
1 &0 &0  \\
0 &p &-1  \\
0 &1 &0
\end{array}
\right)
\end{equation*}
The operation of  removing $\nu (T) $ from $X$ and regluing $T^2\times B^2 $ by the map $\varphi_{p}: S^1 \times T^2 \to \partial \nu (T)$ is called the  $p$  log-transform of $X$ along $T^2$. When $(X, \omega)$ is symplectic and $T^{2}$ is Lagrangian and  $p=\pm1$, this operation preserves symplectic structure and is  equivalent to a {\it Luttinger surgery} up to diffeomorphism (e.g. \cite{a3}). We will refer this operation  as $(a \times b, b)$ Luttinger surgery. Figure 1 describes this as a handlebody operation (cf. \cite{ay})

    \begin{figure}[ht]  \begin{center}  
\includegraphics[width=.72\textwidth]{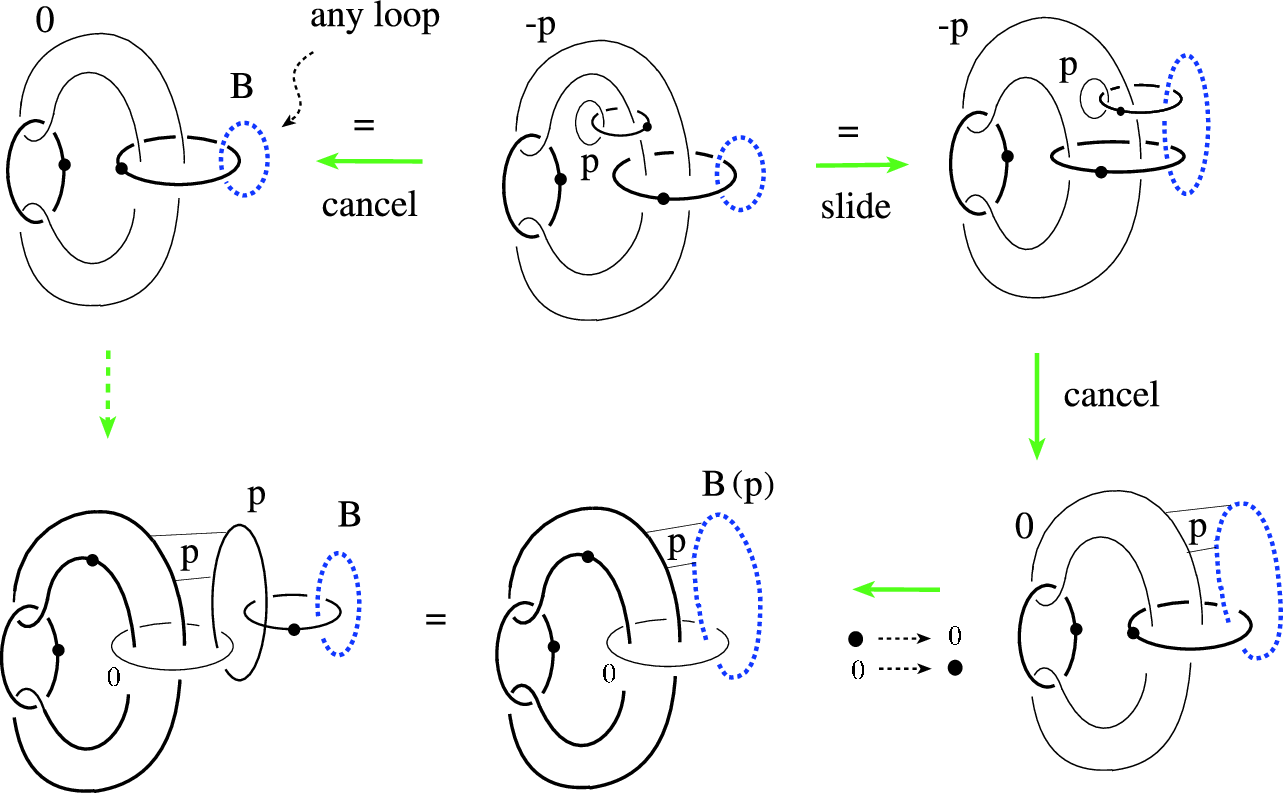}   
\caption{Luttinger surgery ($p=\pm 1$)} 
\end{center}
\end{figure}

\section{Constructing $\tilde{E}_{0}$}

Let $E_{0}=\Sigma_{2}\times T^{2}_{0}$ be the surface of genus two crossed by the punctured torus. Recall that Figure 2 describes   a handlebody picture of $E_{0}$ and the bounday identification  $f: \partial E_{0}\stackrel{\approx}{\longrightarrow} \Sigma_{2}\times S^1$, as shown in \cite{a1}. The knowledge of where the arcs in the figure of $E_{0}$ (top of Figure 1) thrown by the diffeomorphism $f$ is essential to our construction.  By performing the indicated handle slide to $E_{0}$ (indicated by dotted arrow) in Figure 3, we obtain a second equivalent picture of $E_{0}$. By performing Luttinger surgeries to Figure 3  along $(a_1\times C, a_1)$, $ (b_1\times C, b_1)$ we obtain the first picture of Figure 4, and then by  handle slides obtain the second picture. First by an isotopy then a handle slide to Figure 4 we obtain the first and second pictures in Figure 5.  By a further isotopy we obtain the first picture of Figure 6, and then by Luttinger surgeries to $(a_2\times C, C)$, $(a_2\times D, D)$ we obtain the second picture in Figure 6. By introducing canceling handle pairs we express this last picture by a simpler looking first picture of Figure 7. Then by indicated handle slides we obtain the second picture of Figure 7, which is  $\tilde{E_{0}}$.

\subsection{ Calculating $ \pi_{1}(\tilde{E_{0}}) $}

By the indicted handle slide to Figure 7 we obtain Figure 8, which is another picture of $ \tilde{E_{0}}$. In this picture  we also indicated the generators of its fundamental group: $\{a,b, c,d,e,f,g,h,k,p,q\}$. We can read off the relations by tracing the attaching knots of the $2$-handles  (starting at the points indicated by small circles). We get the following relations given by the words:
 
 \vspace{.05in}
 
   $$af, \;ab^{-1}, \; dp, \; cqd^{-1}, \; kh^{-1}, \; befe^{-1}, \; he^{-1}k^{-1}e,$$   
   $$aqa^{-1}d^{-1}\;, ac^{-1}a^{-1}cp^{-1}, \; ag^{-1}a^{-1}h^{-1}g, \; gf^{-1}kg^{-1}k^{-1}, $$$$dk^{-1}b^{-1}pbk, \; dq^{-1}k^{-1}b^{-1}c^{-1}bk, \; cdc^{-1}bfg^{-1}egf^{-1}e^{-1}b^{-1}d^{-1}$$

  \vspace{.05in}
 
 After eliminating $f,\;b,\;p,\;d,\;k $ by using the obvious short words we get: 
 
  \begin{equation*}
 \pi_{1}(\tilde{E_0}) = \left\langle a,\; c,\; e,\; g,\; h,\; q \; \left| 
\begin{array}{cc}
[a,e]=1,  \;  [h,e]=1,  \;  [a,q]=c, & \\
 $$[c^{-1},a]=cq, \; [g,a]=h, \; [g^{-1},h]=a,
  & \\
$$[cq,ah]=1,\; [c,ah]=1,\;  [q^{-1},c][g^{-1},e]=1&
\end{array}  \right\rangle \right.
\end{equation*}

\vspace{.1in}

\section{Constructing $\tilde{E}_{2}$ and $M$}

The first picture of Figure 9 is the handlebody of $S^{3}_{0}(K) \times S^{1}$ (\cite{a2}), where in this picture the horizontal $T_{2}\times D^2$ is clearly visible, but not the vertical torus $T^{2}_{1}$ (which consits of the Seifet surface of $K$ capped off by the $2$-handle given by the zero framed trefoil). In the second picture of Figure 9 we redraw this handlebody so that both vertical and horizontal tori are clearly visible (reader can check this by canceling $1$- and $2$- handle pairs from the second picture to obtain the first picture).  The first picture of Figure 10  gives $(S^{3}_{0}(K) \times S^{1} )\# 2\bar{\C\P}^2$, where by sliding the $2$-handle of $T_{1}$ over the $2$-handle of $T_{2}$ (and by sliding over the two $\bar{{\C\P}}^{1}$'s of the $\bar{{\C\P}}^{2}$'s) we obtained the imbedding  $\Sigma_{2}\times D^2 \subset (S^{3}_{0}(K) \times S^{1} )\# 2\bar{\C\P}^2$. By isotopies and handle slides (indicated by dotted arrows) we obtain the second picture in Figure 10, and then the Figures 11 and 12. Either  pictures of Figure 12 represent handlebodies of $(S^{3}_{0}(K) \times S^{1} )\# 2\bar{\C\P}^2$ (both have different advantages).  Figure 13 is the same as the first picture of Figure 12, drawn in an exaggerated way so that $\Sigma_{2} \times D^2$ is clearly visible. We now want to remove this  $\Sigma_{2} \times D^2$  from inside this handlebody of $(S^{3}_{0}(K) \times S^{1} )\# 2\bar{\C\P}^2$ and replace it with $\tilde{E_{0}}$. The arcs in Figure 2 (describing the diffeomorphism $f$) and also in all the subsequent Figures 3 to 7 show us how to do this, resulting with the handlebody picture of $M$ in Figure 14. Clearly Figure 15 is another handlebody of $M$ (where we used the second picture of Figure 12 instead of the first). 

\subsection{ Checking $ \pi_{1}(M)=0 $}

Clearly we can calculate  $ \pi_{1}(M)$ from the presentation of $ \pi_{1}(\tilde{E_{0}}) $ from Figure 8 (Section 2.1) by introducing new generators $x$ and $y$ (Figure 13 and 14) and $7$ new relations coming from the new $2$-handles in Figure 14. The $7$ new relations are given by the words:
$$(xyx)^{-1}(yxy),\; q^{-1}y,\;  g^{-1}xy^{-1}x^{-1}y,$$$$ e^{-1}yx^{-1},\; cyx^{-1}c^{-1}xy^{-1},\; ege^{-1}g^{-1},\; a^{-1}h^{-1}g^{-1}hg$$

 After eliminating $x$ and $y$ from the two short relations we get:  
 
   \begin{equation*}
 \pi_{1}(M) = \left\langle a,\; c,\; e,\; g,\; h,\; q \; \left| 
\begin{array}{cc}
[a,e]=1,  \;  [h,e]=1,  \;  [a,q]=c, & \\
 $$[c^{-1},a]=cq, \; [g,a]=h, \; [g^{-1},h]=a,
  & \\
$$[cq,ah]=1,\; [c,ah]=1,\;  [q^{-1},c][g^{-1},e]=1&\\
a=[h^{-1}g^{-1}], \; [e,g]=1,\; [c,e]=1&\\
g=e^{-1}qeq,\; q^{2}=eqe^{-1}qe\\
\end{array}  \right\rangle \right.
\end{equation*}
From this one can easily check $ \pi_{1}(M) =0$. For example, since $q=c^{-1}[c^{-1},a]$, and $a$ and $c$ commutes with $e$, then $q$ commutes with $e$. Hence the last two relations imply $e=1$ and $g=q^2$.  The relations $[cq,ah]=1$ and $[c,ah]=1$ imply $[q,ah]=1$, hence $[g,ah]=1$, in tern this together with $ [g^{-1},h]=a$ gives $a=1$, then relations $[a,q]=c$  and $[g,a]=h$ gives $c=1$ and $h=1$, and hence $[c^{-1},a]=cq$  implieq $q=1$ and so $g=q^2=1$.

    \begin{figure}[ht]  \begin{center}  
\includegraphics[width=.9\textwidth]{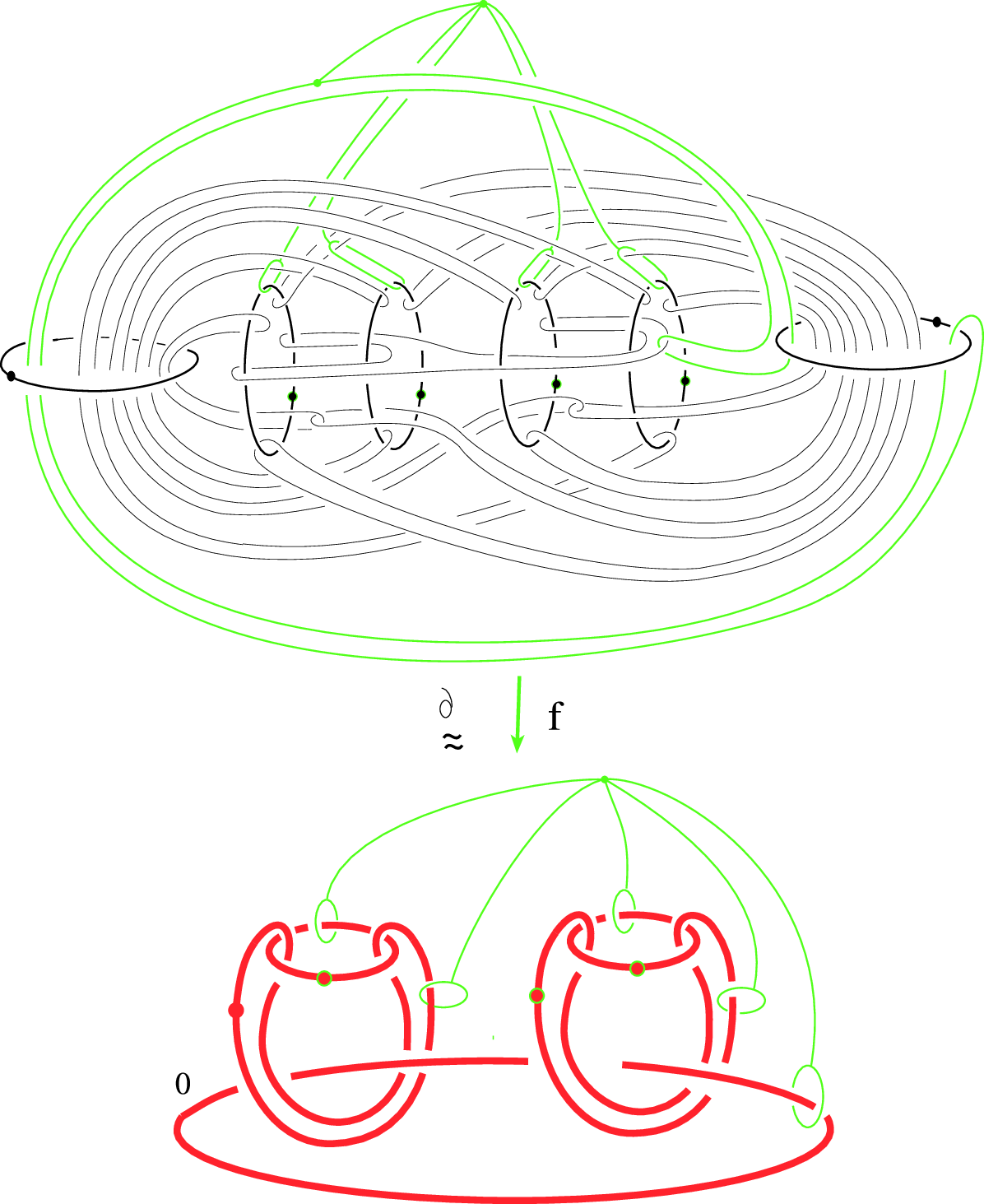}   
\caption{$E_{0}$ and  $f: \partial E_{0}\stackrel{\approx}{\longrightarrow} \Sigma_{2}\times S^1$} 
\end{center}
\end{figure} 

    \begin{figure}[ht]  \begin{center}  
\includegraphics[width=.9\textwidth]{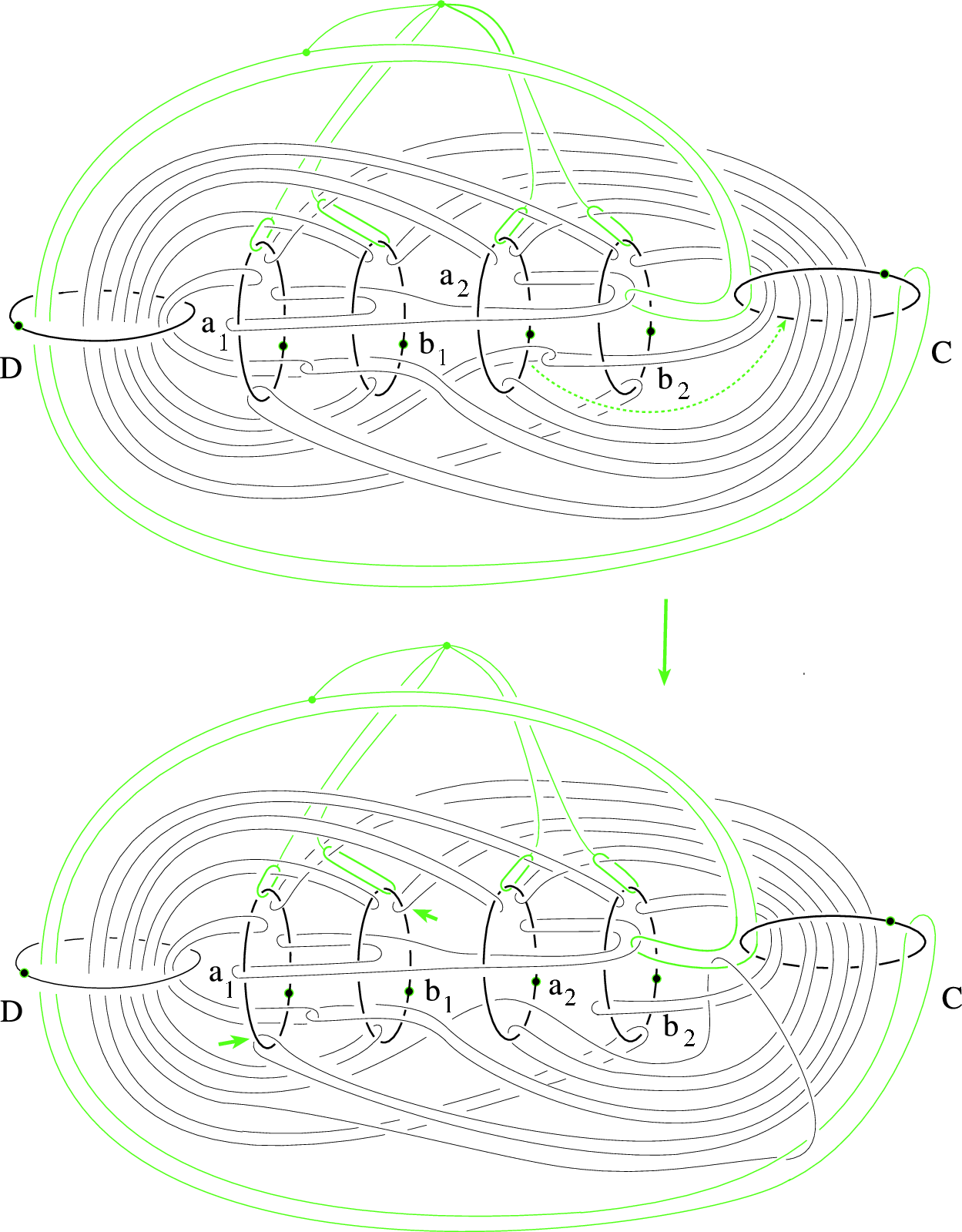}   
\caption{Handle slide (the pair of thick  arrows  indicate where we will perform Luttinger surgeries next)} 
\end{center}
\end{figure}

   \begin{figure}[ht]  \begin{center}  
\includegraphics[width=.9\textwidth]{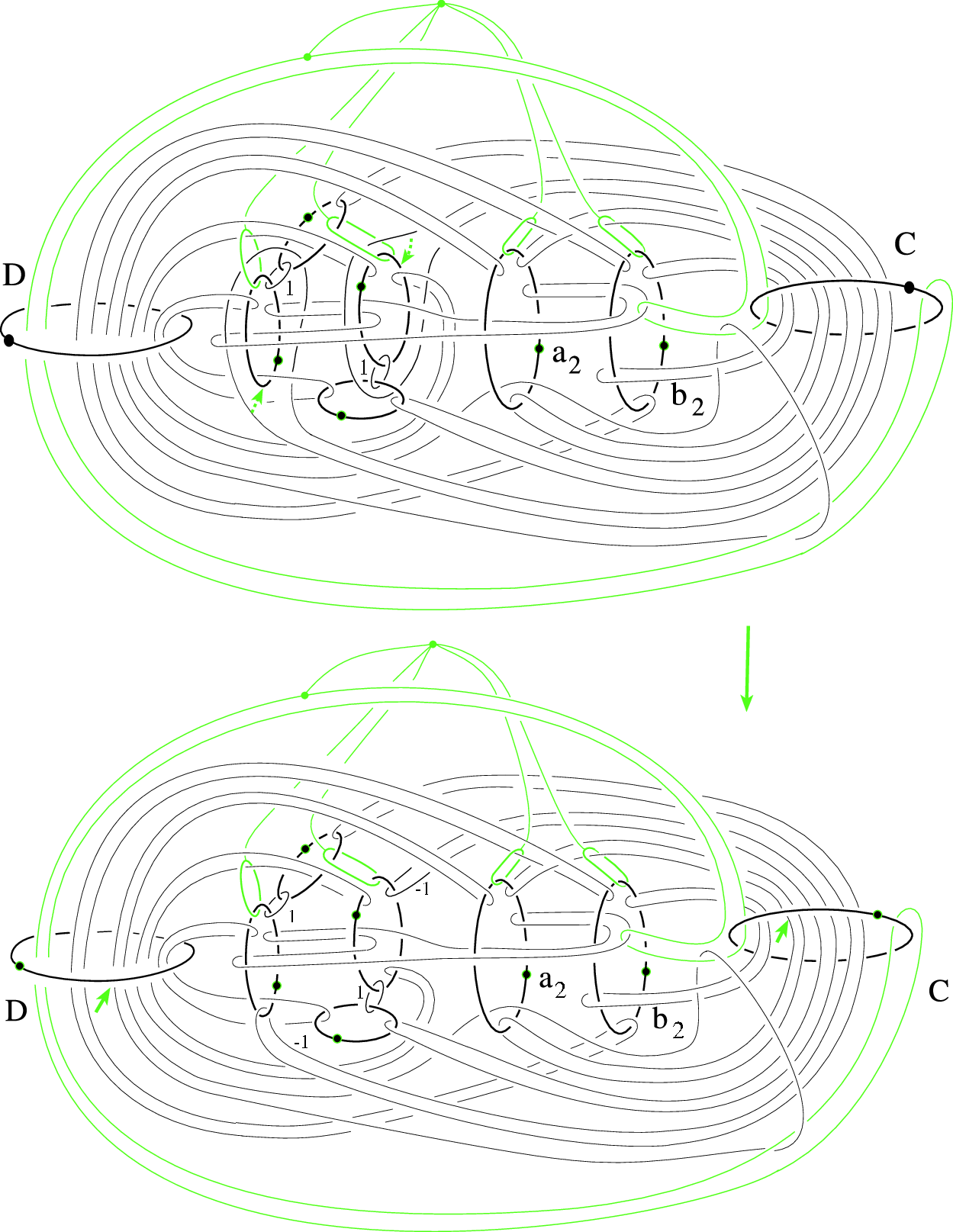}   
\caption{First we performed Luttinger surgeries along $( a_1\times C, a_1)$, $( b_1\times C, b_1)$, then handle slides (thick arrows indicate where we will perform Luttinger surgeries next)} 
\end{center}
\end{figure} 

   \begin{figure}[ht]  \begin{center}  
\includegraphics[width=.9\textwidth]{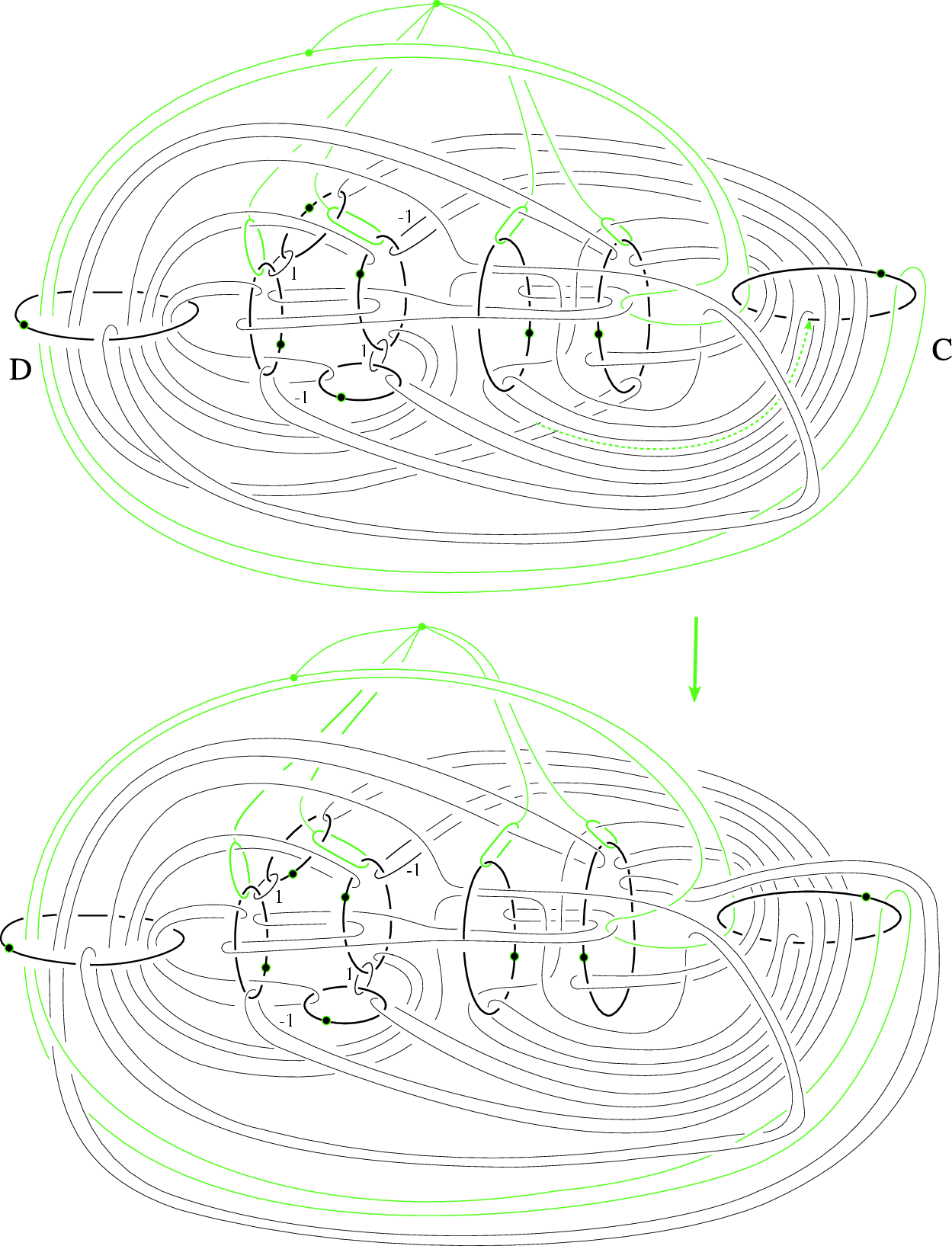}   
\caption{Isotopy and a handle slides} 
\end{center}
\end{figure}

   \begin{figure}[ht]  \begin{center}  
\includegraphics[width=.9\textwidth]{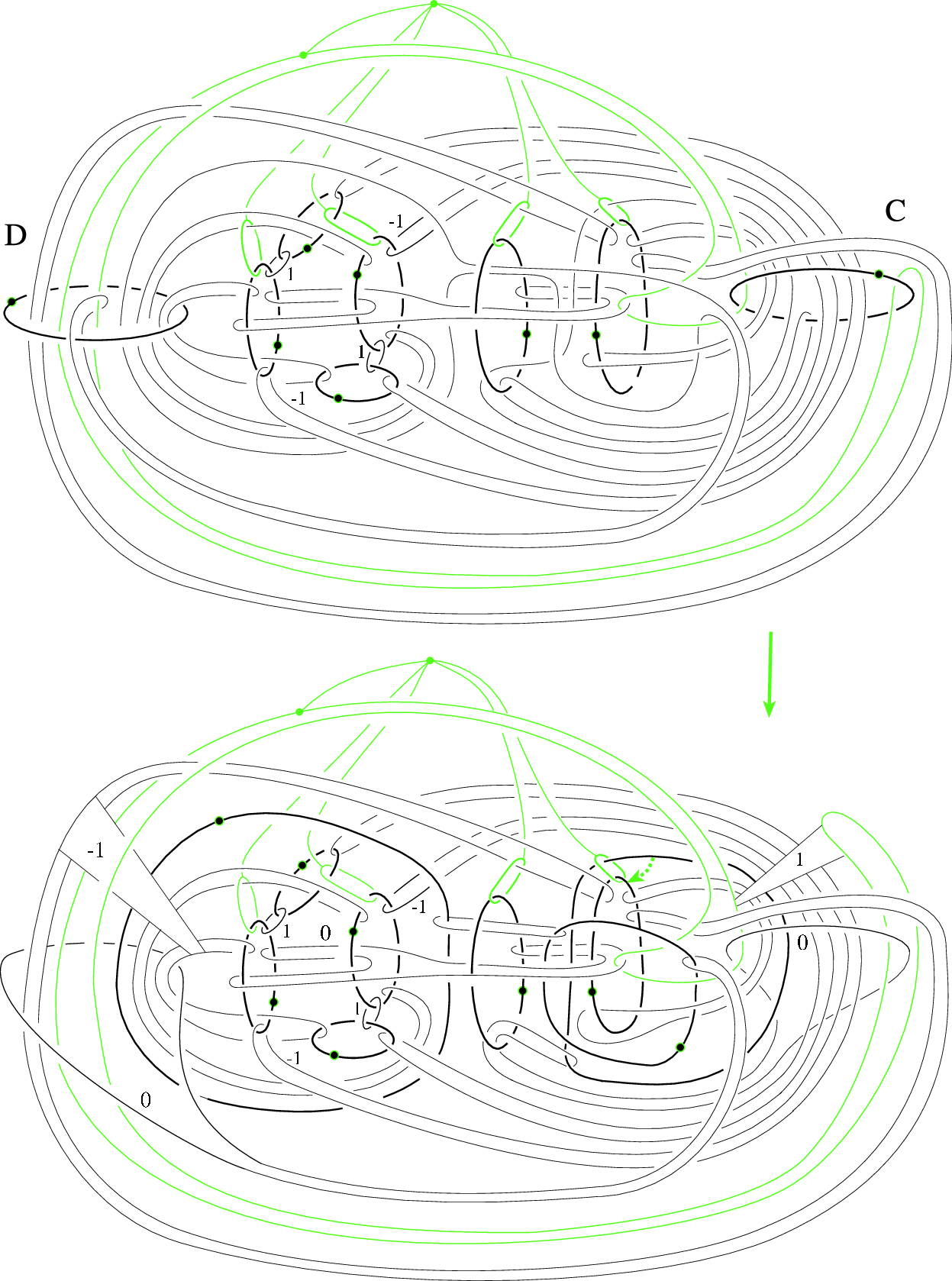}   
\caption{Luttinger surgeries along $(a_2 \times C, C)$, $( a_2\times D, D)$ } 
\end{center}
\end{figure}

   \begin{figure}[ht]  \begin{center}  
\includegraphics[width=.95\textwidth]{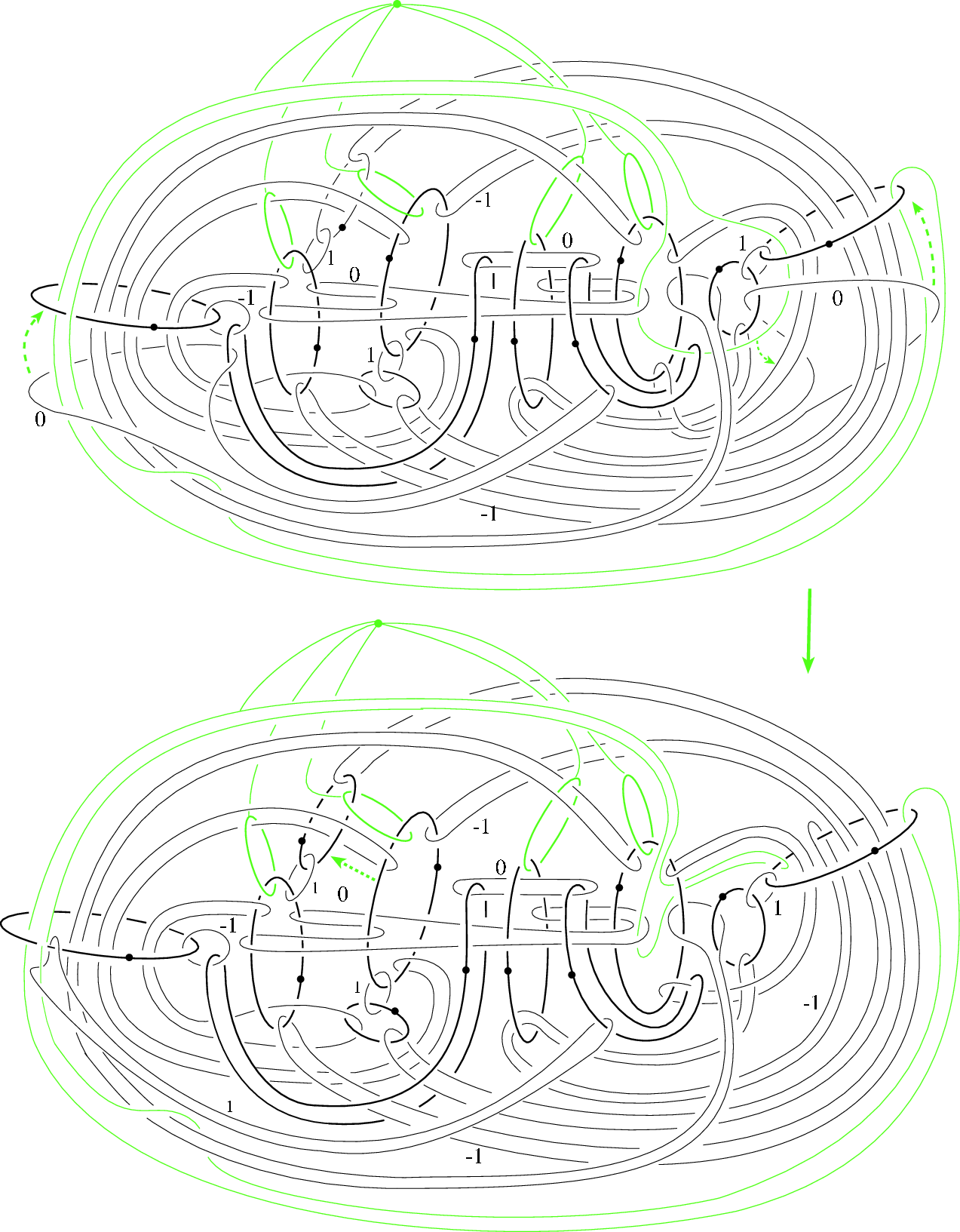}   
\caption{More isotopies and handle slides and getting $\tilde{E_{0}}$}
\end{center}
\end{figure} 

   \begin{figure}[ht]  \begin{center}  
\includegraphics[width=1\textwidth]{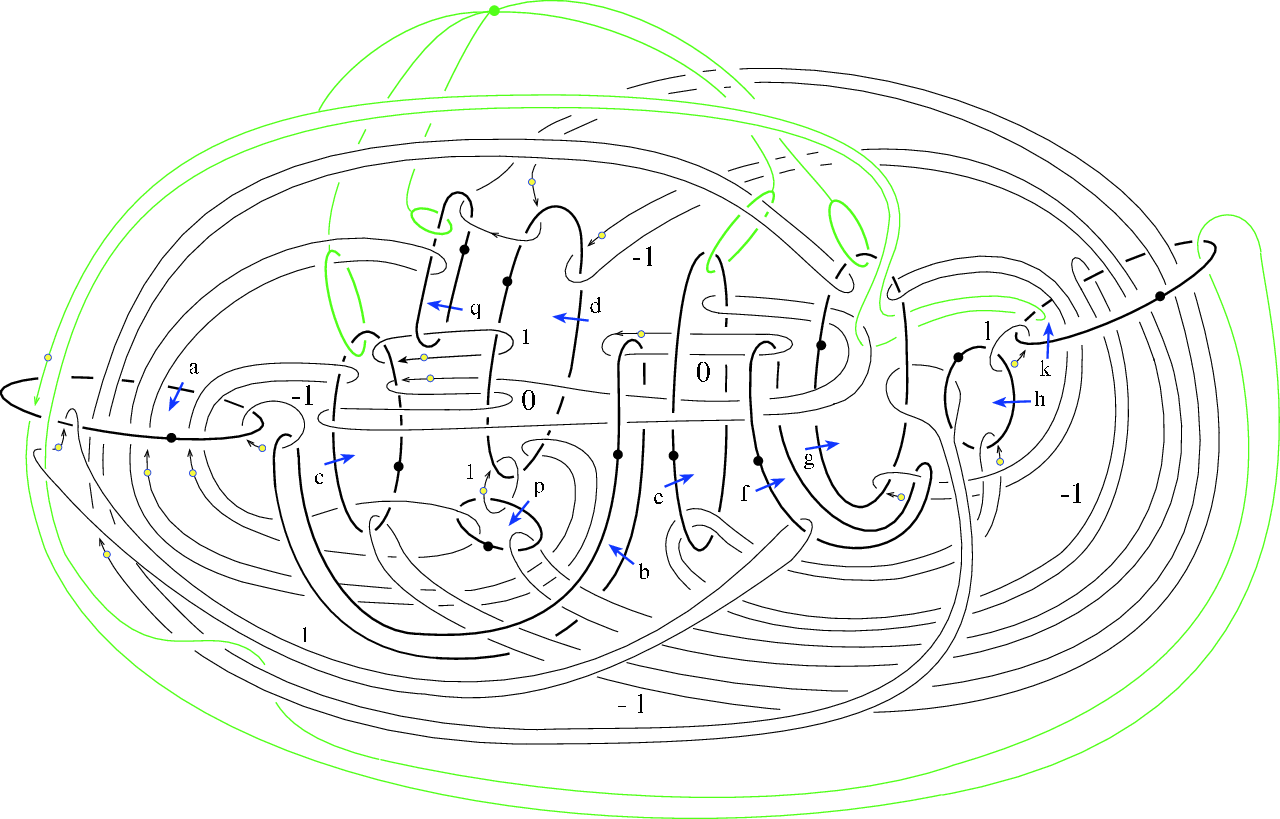}   
\caption{$\tilde{E_{0}}$, with generators of $\pi_{1} (\tilde{E_{0}})$ are indicated } 
\end{center}
\end{figure} 

   \begin{figure}[ht]  \begin{center}  
\includegraphics[width=.7\textwidth]{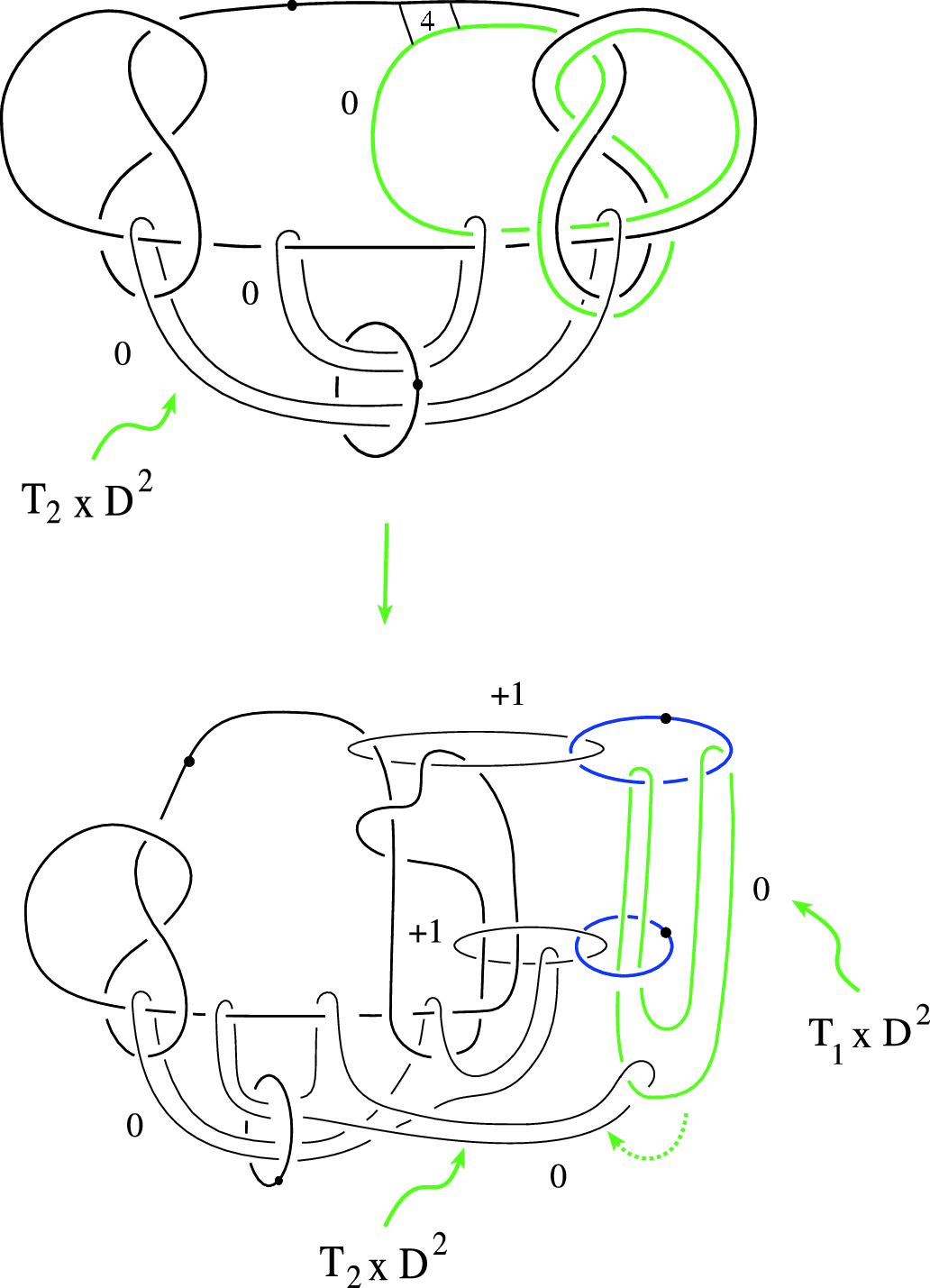}   
\caption{$S^{3}_{0}(K) \times S^{1}$}
\end{center}
\end{figure} 

   \begin{figure}[ht]  \begin{center}  
\includegraphics[width=.7\textwidth]{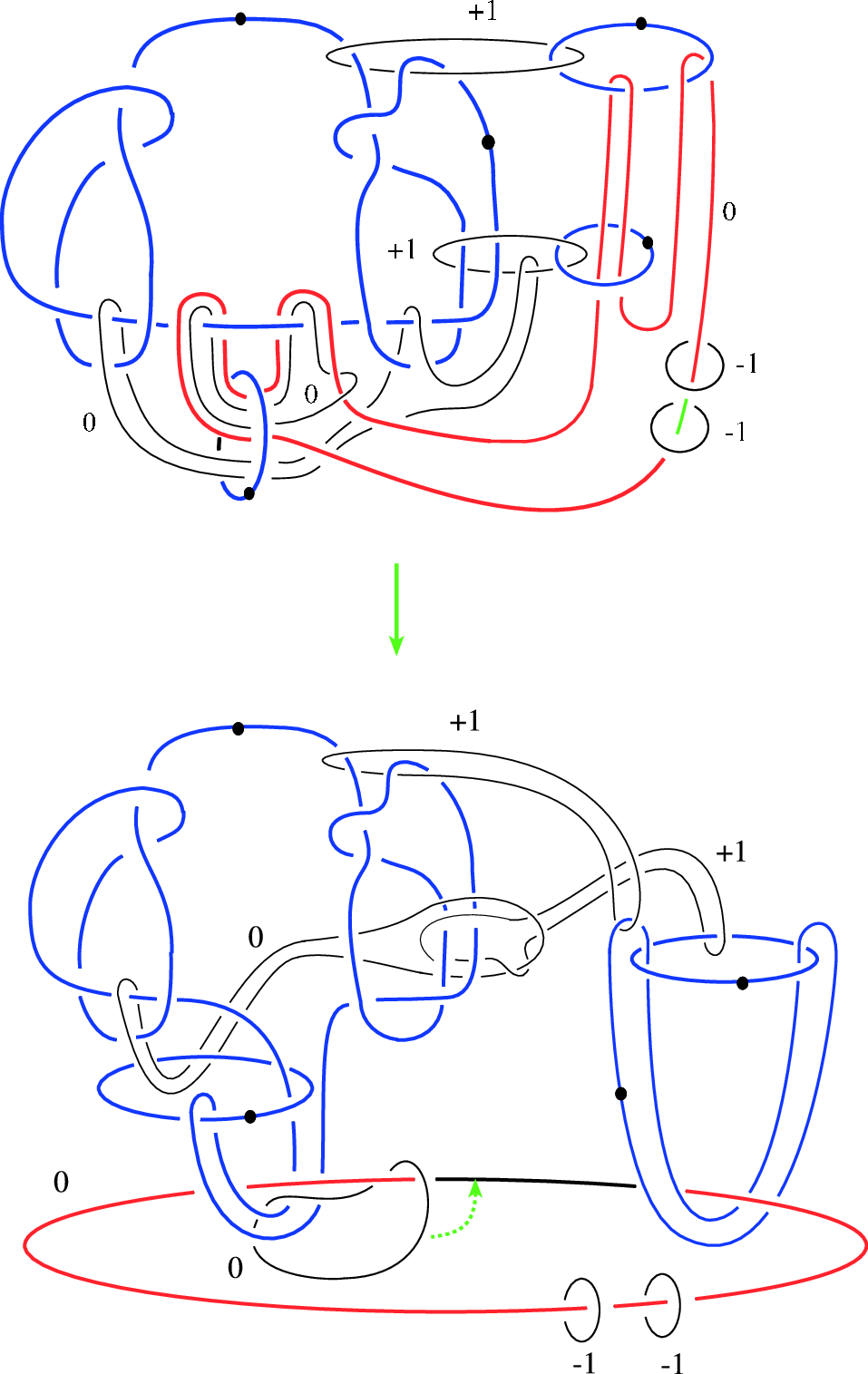}   
\caption{$(S^{3}_{0}(K) \times S^{1} )\# 2\bar{\C\P}^2$} 
\end{center}
\end{figure} 

   \begin{figure}[ht]  \begin{center}  
\includegraphics[width=.7\textwidth]{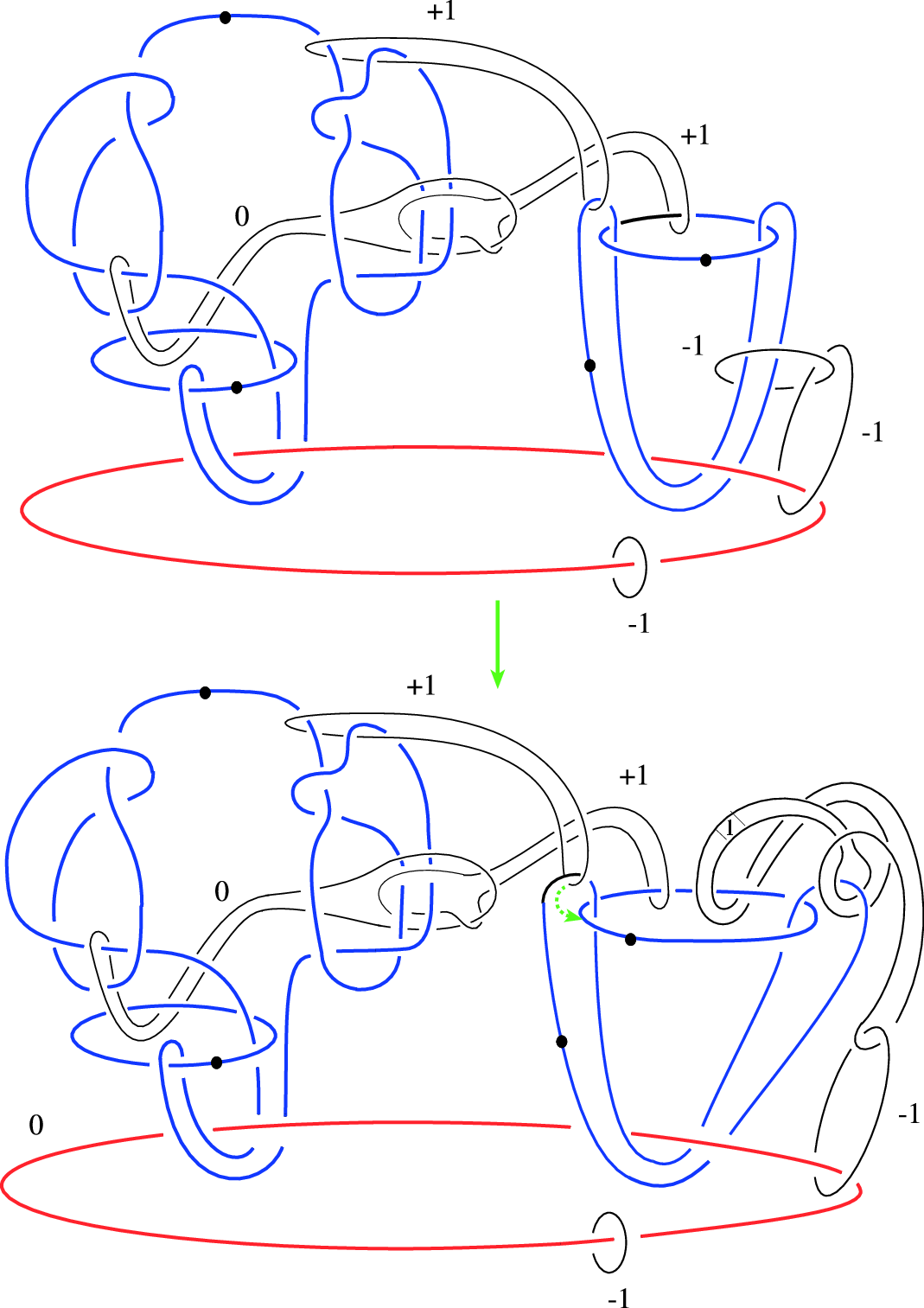}   
\caption{$(S^{3}_{0}(K) \times S^{1} )\# 2\bar{\C\P}^2$} 
\end{center}
\end{figure} 

   \begin{figure}[ht]  \begin{center}  
\includegraphics[width=.7\textwidth]{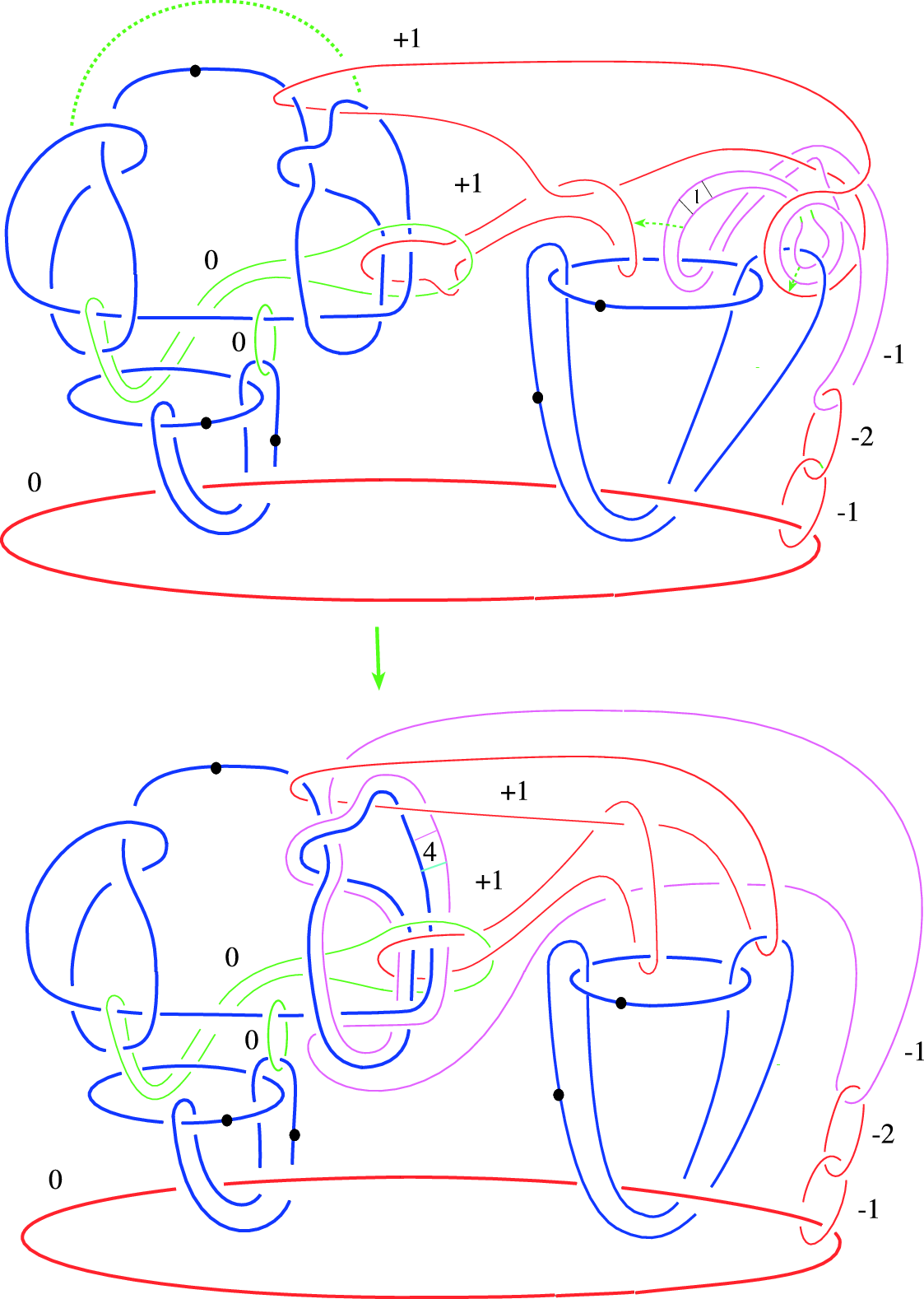}   
\caption{$(S^{3}_{0}(K) \times S^{1} )\# 2\bar{\C\P}^2$} 
\end{center}
\end{figure} 

   \begin{figure}[ht]  \begin{center}  
\includegraphics[width=.7\textwidth]{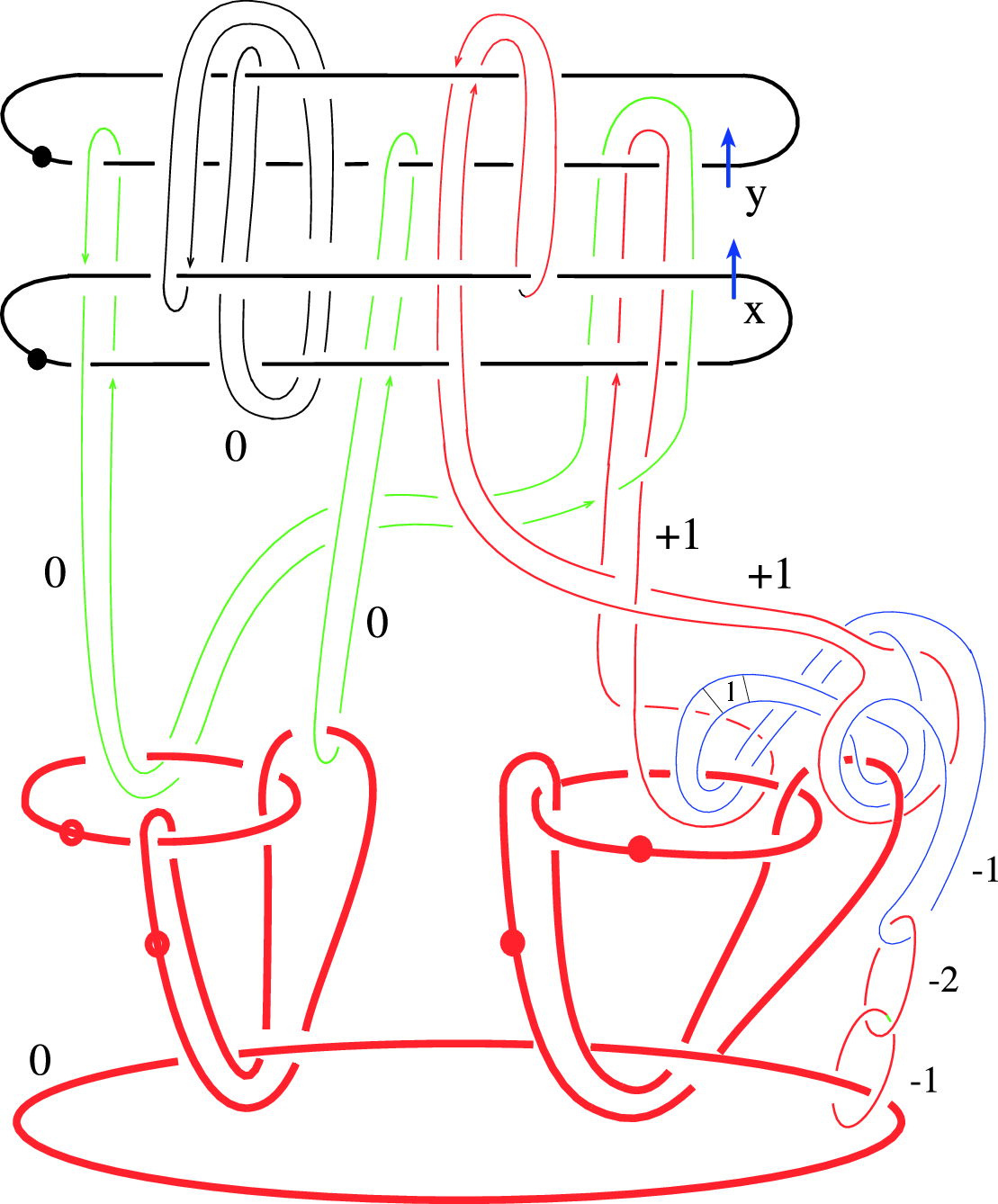}   
\caption{$(S^{3}_{0}(K) \times S^{1} )\# 2\bar{\C\P}^2$} 
\end{center}
\end{figure} 

   \begin{figure}[ht]  \begin{center}  
\includegraphics[width=.99\textwidth]{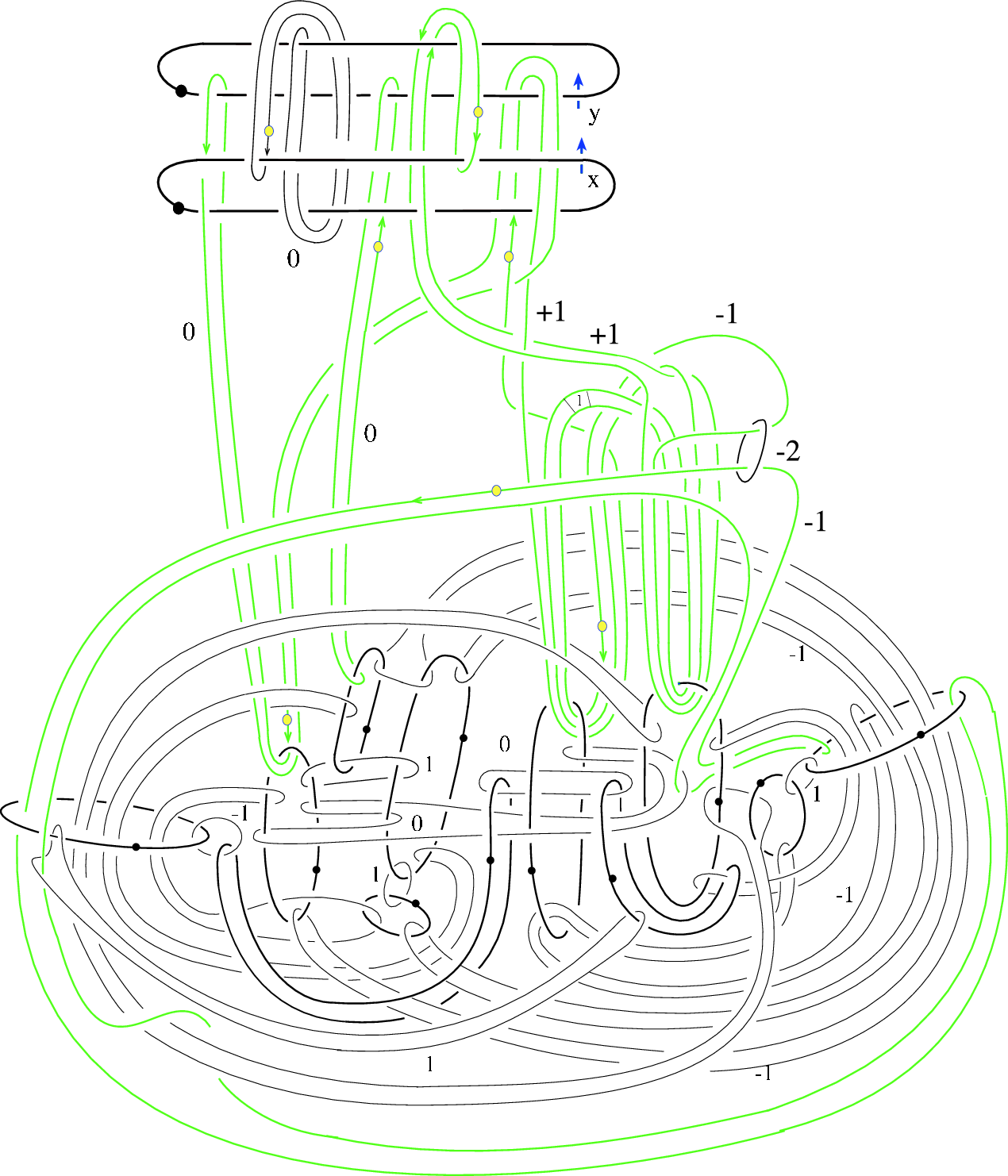}   
\caption{ $M$} 
\end{center}
\end{figure}

   \begin{figure}[ht]  \begin{center}  
\includegraphics[width=.99\textwidth]{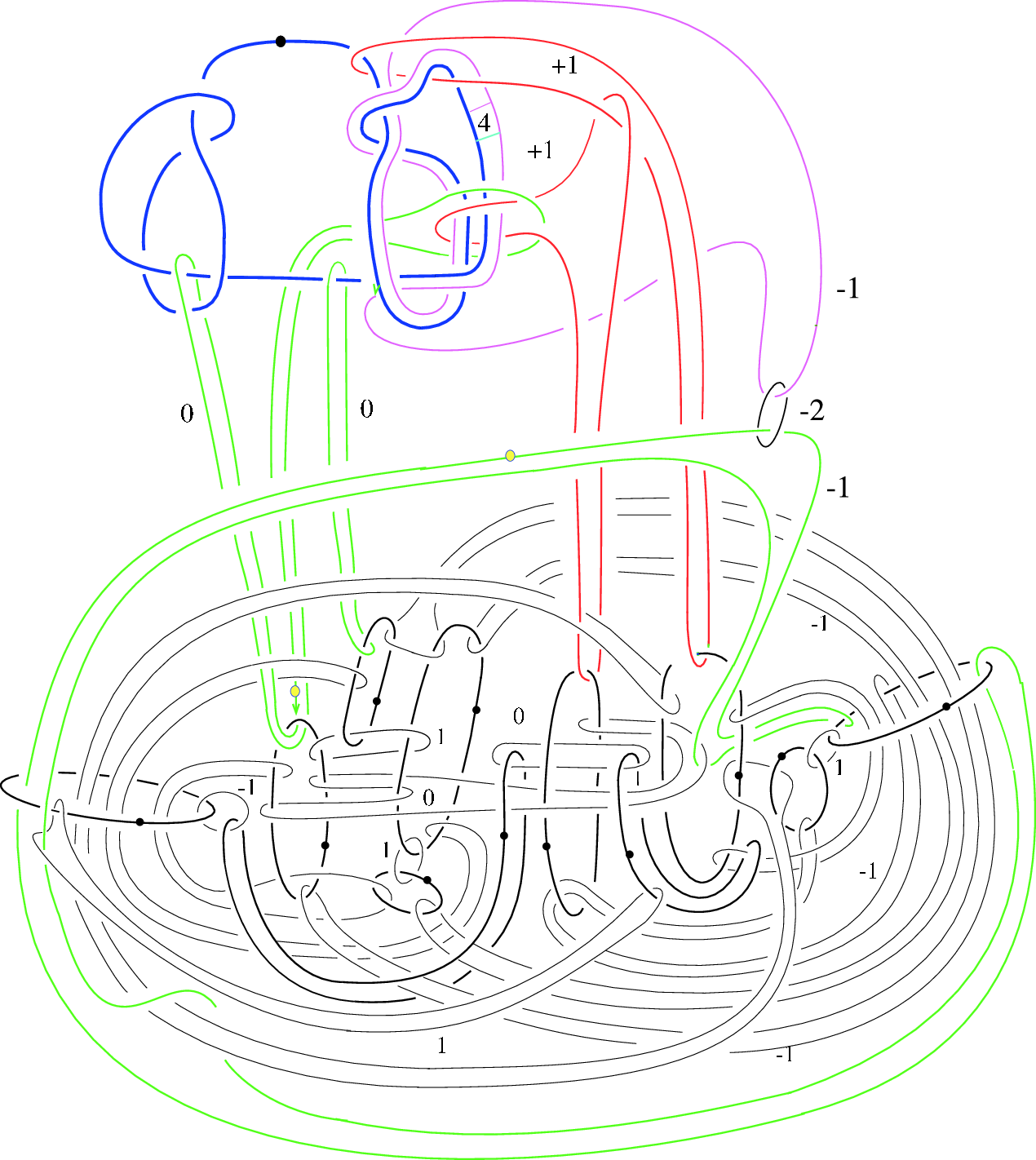}   
\caption{ $M$} 
\end{center}
\end{figure}


\clearpage 

\begin{thebibliography}{9999}

\bibitem  [A1] {a1} S. Akbulut, {\em Catanese-Ciliberto-Mendes Lopes surface,}
	arXiv:1101.3036.
	
	 \bibitem [A2]{a2}  S. Akbulut, {\em  A fake cusp and a fishtail,}   Turkish Jour. of Math 1 (1999), 19-31.   arXiv:math.GT/9904058.
	 
	 	  \bibitem [A3]{a3}  S. Akbulut, {\em 4-Manifolds}, Book in preparation, available from\\
  \url{http://www.math.msu.edu/~akbulut/papers/akbulut.lec.pdf}, 2014.

	 	 	 
	\bibitem [Ak]{ak}A. Akhmedov  {\em Small exotic 4-manifolds}, Algebraic and Geometric Topology, 8 (2008), 1781-1794.
	 
	  \bibitem [AP]{ap} A. Akhmedov, B. D. Park, {\em Exotic 4-manifolds with small Euler characteristic}, Inventiones Mathematicae,
173 (2008), 209-223.

 \bibitem[AY]{ay} S. Akbulut and K. Yasui, {\em Corks, Plugs and exotic structures}, \\Journal of G\"{o}kova Geometry Topology, volume \textbf{2} (2008), 40--82.  

	  
	\bibitem  [FPS]{fps} R. Fintushel, B. D. Park, R. J. Stern {\em Reverse engineering small 4-manifolds}, Algebraic and Geometric
Topology, 8 (2007), 2103-2116.

\bibitem[GS]{gs} R. E. Gompf and A.I. Stipsicz,  {\em $4$-manifolds and Kirby calculus}, Graduate Studies in Mathematics 20. American Mathematical Society, 1999.
	  
\bibitem [BK]{bk} S.  Baldridge, P. Kirk  {\em A symplectic manifold homeomorphic but not diffeomorphic to ${\C\P}^{2}\# 3 {\C\P}^{2}$}, Geometry
and Topology 12, (2008), 919-940.


\end{thebibliography}
\end{document}